\documentclass{article}

\usepackage{multicol} 
\usepackage[svgnames]{xcolor} 

\usepackage{times} 

\usepackage{graphicx} 
\graphicspath{{figures/}} 
\usepackage{booktabs} 
\usepackage[font=small,labelfont=bf]{caption} 
\usepackage{amsfonts, amsmath, amsthm, amssymb} 
\usepackage{wrapfig} 

\title{\bf Generalized Fibonacci Sequences and Binet-Fibonacci Curves}

\author{Merve \"{O}zvatan and Oktay K. Pashaev\\ \\
Department of Mathematics \\
Izmir Institute of Technology \\
 Izmir, 35430, Turkey}

\begin{document}

\maketitle

\begin{abstract}
We have studied several generalizations of Fibonacci sequences as the sequences with arbitrary initial values, the addition of two and more Fibonacci subsequences and Fibonacci polynomials with arbitrary bases. For Fibonacci numbers with congruent indices we derived general formula in terms of generalized Fibonacci polynomials and Lucas numbers. By extending Binet formula to arbitrary real numbers we constructed Binet-Fibonacci curve in complex plane. For positive numbers the curve is oscillating with exponentially vanishing amplitude, while for negative numbers it becomes Binet Fibonacci spiral. Comparison with the nautilus curve shows quite similar behavior. Areas under the curve and curvature characteristics are calculated, as well as the asymptotic of relative characteristics. Asymptotically for n going to infinity the region become, infinitely wide with infinitesimally small height so that the area has finite value $A_{\infty}=\frac{2}{5 \pi} ln(\varphi)$.

\end{abstract}





\section{Introduction}
\subsection{Fibonacci Sequence}
\quad Fibonacci sequence is defined by recursion formula;
\begin{equation}
F_{n+1}=F_{n}+F_{n-1} \label{recursion},
\end{equation}
 where $F_0=0$  ,  $F_1=1$ , $n=1,2,3,...$
 First few Fibonacci numbers are : $$ 0,1,1,2,3,5,8,13,21,34,55,89,144,233,377,610,987...  $$

 The sequence is named after Leonardo Fibonacci(1170-1250) \cite{fibonacci}.
 Fibonacci numbers appear in Nature so frequently that they can be considered as Nature's Perfect Numbers.
 Also, another important Nature's number,  the  Golden ratio, which seen, in every area of life and art, and usually it is associated with aesthetics, is related to Fibonacci sequence.

\subsection{Binet Formula}
\quad Formula for Fibonacci sequence was derived by Binet in 1843. It has the form;
\begin{equation}
F_n=\frac{\varphi^n-\varphi'^n}{\varphi-\varphi'} \label{binetformula}
\end{equation}
where $\varphi$ and $\varphi'$ are roots of equation $x^2-x-1=0$, having the values;\\

$$\varphi=\frac{1+\sqrt{5}}{2}\approx 1,6180339.. \,\,\,\, \mbox{and} \,\,\,\, \varphi'=\frac{1-\sqrt{5}}{2}\approx -0,6180339..$$\\

From this formula, due to irrational character of Golden ratio $\varphi$ and $\varphi'$, it is not at all evident that $F_n$ are integer numbers. Though it is clear from the recursion formula (\ref{recursion}).\\

Number $\varphi$ is called the Golden ratio. Binet formula allows one to find corresponding Fibonacci numbers directly without of using recursion, like $F_{98}$ ,$F_{50}$ ,... For example, to find $F_{20}$ by using Binet formula we have; \\
\begin{equation}
F_{20}=\frac{\varphi^{20}-\varphi'^{20}}{\varphi-\varphi'}=6765. \nonumber \\
\end{equation}

Binet formula allows one to define also Fibonacci numbers for negative n;\\

\begin{equation}
F_{-n}=\frac{\varphi^{-n}-\varphi'^{-n}}{\varphi-\varphi'}=\frac{\frac{1}{\varphi^n}-\frac{1}{\varphi'^n}}{\varphi-\varphi'}=\frac{\varphi'^{n}-\varphi^{n}}{\varphi-\varphi'}\frac{1}{(\varphi\varphi')^n} \:,\nonumber\\
\end{equation}\\

since $\varphi\varphi'=-1$,\\
\begin{equation}
\Rightarrow \qquad=-\frac{\varphi^{n}-\varphi'^{n}}{\varphi-\varphi'}\frac{1}{(-1)^n}=-F_{n}\frac{1}{(-1)^n}=-F_{n}(-1)^n=(-1)^{n+1}F_{n}  \nonumber \\
\end{equation}
So, we have;\\
\begin{equation}
F_{-n}=(-1)^{n+1}F_{n} \label{relationwith negative index}
\end{equation}

\section{Generalized Fibonacci numbers}
\qquad Here, we are going to study different generalizations of Fibonacci numbers.
 As a first generalization, by choosing different initial numbers $G_0$ and $G_1$, but preserving the recursion formula (\ref{recursion}), we can define so called generalized Fibonacci numbers. For example, if $G_0=0$, $G_1=4$ we have the sequence;\\
$$4,4,8,12,20,32,52,...$$\\\qquad
\subsection{Addition of two Fibonacci sequences}
Let us we consider generalized Fibonacci numbers, with the recursion formula $G_{n+1}=G_{n}+G_{n-1}$, and an arbitrary initial numbers $G_0$ $\&$ $G_1$;\\
$G_0=G_0$ \\
$G_1=G_1$ \\
$G_2=G_1+G_0=F_2G_1+F_1G_0$ \\
$G_3=G_2+G_1=2G_1+G_0=F_3G_1+F_2G_0$ \\
$G_4=G_3+G_2=3G_1+2G_0=F_4G_1+F_3G_0$ \\
$G_5=G_4+G_3=5G_1+3G_0=F_5G_1+F_4G_0$ \\
\phantom{abc}.\\
\phantom{abc}.\\
Then
\begin{equation}
G_n=G_1F_n+G_0F_{n-1} \label{2rec.for.}
\end{equation}
 is obtained as a recursion relation. This shows that $G_n$ is a linear combination of two Fibonacci sequences. \\

By substituting Binet formulas for $F_n,F_{n-1}$ gives;\\
\begin{align}
G_n&=G_1F_n+G_0F_{n-1}, \nonumber \\
G_n&=G_1\frac{\varphi^n-\varphi'^n}{\varphi-\varphi'}+G_0\frac{\varphi^{n-1}-\varphi'^{n-1}}{\varphi-\varphi'}, \nonumber \\
G_n&=\frac{1}{\varphi-\varphi'}\left[(G_1\varphi^n+G_0\varphi^{n-1})-(G_1\varphi'^n+G_0\varphi'^{n-1})\right], \nonumber \\
G_n&=\frac{1}{\varphi-\varphi'}\left[(G_1\varphi^n+G_0\varphi^{n}\frac{1}{\varphi})-(G_1\varphi'^n+G_0\varphi'^{n}\frac{1}{\varphi'})\right],\,\, \mbox{and since\,\,} \varphi'=-\frac{1}{\varphi}, \nonumber \\
G_n&=\frac{1}{\varphi-\varphi'}\left[(G_1\varphi^n-G_0\varphi^{n}\varphi')-(G_1\varphi'^n-G_0\varphi'^{n}\varphi)\right] \nonumber \\
G_n&=\frac{(G_1-\varphi'G_0)\varphi^n-(G_1-\varphi G_0)\varphi'^n}{\varphi-\varphi'}. \label{bin.form.for.gen.fib.}
\end{align} \\
\phantom{abh}Equation (\ref{bin.form.for.gen.fib.}) we called Binet type formula for generalized Fibonacci numbers. Also note that, if $G_0=0$ and $G_1=1$, our equation becomes Binet formula (\ref{binetformula}) for Fibonacci numbers.\\

Now, let us see in other way that our recursion  $G_{n+1}=G_{n}+G_{n-1}$ is valid for the equation (\ref{bin.form.for.gen.fib.}). Suppose;\\
$$G_{n+1}=AG_{n}+BG_{n-1}.$$\\
\phantom{abh}We have to find A=B=1. By writing corresponding recursion relation formulas for generalized Fibonacci numbers; \\
$$G_1F_{n+1}+G_0F_{n}=A(G_1F_n+G_0F_{n-1})+B(G_1F_{n-1}+G_0F_{n-2}),$$\\
$$G_1(F_n+F_{n-1})+G_0(F_{n-1}+F_{n-2})=G_1(AF_{n}+BF_{n-1})+G_0(AF_{n-1}+BF_{n-2}).$$
From above equality, we can say that A and B are 1.\\

\subsection{Arbitrary linear combination}

As we have seen, the linear combination of two Fibonacci sequences (\ref{2rec.for.}) with shifted by 1 index $n$, produces Generalized Fibonacci numbers. Now let us think more general case of recursion formula for Generalized Fibonacci numbers defined as; $$G_n^{(k)}=\alpha_0F_n+\alpha_1F_{n+1}+...+\alpha_kF_{n+k}.$$\\
Here the sequence $G_n^{(k)}$ is determined by coefficient, of a polynomial degree k; \\
$$P_k(x)=\alpha_0+\alpha_1 x+...+\alpha_k x^{k}.$$\\
As easy to see, the sequence $G_{n+1}^{(k)}$ satisfies the standard recursion formula \\
\begin{equation}
G_{n+1}^{(k)}=G_n^{(k)}+G_{n-1}^{(k)}.
\end{equation}  \\
By using Binet Formula for Fibonacci numbers we can derive Binet type formula for our Generalized Fibonacci numbers as follows;\\

\begin{equation}
G_n^{(k)}=\frac{\varphi^nP_k(\varphi)-\varphi'^nP_k(\varphi')}{\varphi-\varphi'} \label{implct.form.of.most.gen.case.for.shfted.gen.nbs}.\\
\end{equation}

\subsubsection{Binet formula and dual polynomials}
The above Binet type formula (\ref{implct.form.of.most.gen.case.for.shfted.gen.nbs}) can be rewritten in terms of dual polynomials.\\
\textbf{Definition:}\quad For given polynomial degree k,
\begin{equation}
P_k(x)=\sum^k_{l=0}a_l x^l
\end{equation}
we define the dual polynomial degree k according to formula
\begin{equation}
\tilde P_k(x)=\sum^k_{l=0}a_{k-l} x^l.
\end{equation}
It means that if $P_k(x)$ is a vector with components $(a_0,a_1,...,a_k)$, then the dual $\tilde P_k(x)$ is the vector with components $(a_k,a_{k-1},...,a_0)$.
These polynomials are related by formula;
\begin{equation}
P_k(x)=x^k \tilde P_k\left(\frac{1}{x}\right).
\end{equation}

Let us consider polynomial:\\

$$P_k(\varphi)=\alpha_0+\alpha_1 \varphi+...+\alpha_k \varphi^{k},$$ \\
then
$$P_k(\varphi)=\varphi^k\left(\alpha_0\frac{1}{\varphi^k}+\alpha_1\frac{1}{\varphi^{k-1}}+...+\alpha_k\right)\;\;\;\mbox{or by using}\;\;\;\varphi'=-\frac{1}{\varphi};$$ \\
$$P_k(\varphi)=\varphi^k\left(\alpha_0(-\varphi')^{k}+\alpha_1(-\varphi')^{k-1}+...+\alpha_k\right),$$ \\
and
$$P_k(\varphi)=\varphi^k\left(\alpha_k+\alpha_{k-1}(-1)\varphi'+\alpha_{k-2}(-1)^2\varphi'^2+...+\alpha_1(-1)^{k-1}\varphi'^{k-1}+   \alpha_0(-1)^{k}\varphi'^k\right)$$ \\
Thus, finally we have;\\

$$P_{\alpha_0,\alpha_1,...,\alpha_k}(\varphi)=\varphi^k \tilde P_{\alpha_k,-\alpha_{k-1},(-1)^2\alpha_{k-2},...,(-1)^{k-1}\alpha_1,(-1)^k\alpha_0}(\varphi').$$\\
It can be written as, \\
$$P_k(\varphi)=\varphi^k \tilde{P}_k(-\varphi'),$$ where $\tilde{P}_k(-\varphi')$ is the dual polynomial to $P_k(\varphi)$. \\
\phantom{ab}By following the same procedure and starting from $P_k(\varphi')$, we can easily obtain that;\\
$$P_k(\varphi')=\varphi'^k \tilde{P}_k(-\varphi).$$\\

Therefore, for given polynomials with arguments $\varphi$ or $\varphi'$ we have the dual one with argument $-\varphi'$ and $-\varphi$ respectively. This allow us to write Binet type formula for Generalized Fibonacci numbers in two different forms:\\

$$G_n^{(k)}=\frac{\varphi^nP_k(\varphi)-\varphi'^nP_k(\varphi')}{\varphi-\varphi'}=\frac{\varphi^{n+k}\tilde{P}_k(-\varphi')-\varphi'^{n+k}\tilde{P}_k(-\varphi)}{\varphi-\varphi'}.$$

\section{Fibonacci Polynomials}
\qquad In previous section we studied Generalized Fibonacci numbers with Fibonacci recursion formula. Here we are going to generalize this recursion formula by introducing two arbitrary numbers. Then, the corresponding sequence will depend on two numbers. This sequence of two variable polynomials is called the Fibonacci polynomials.

Let $\mbox{a}\, \&\, \mbox{b}$ be (real) roots of the second order equation;\\
$(x-a)(x-b) = x^2-(a+b)x+ab=0 $. And let us say $ a+b=p $ and $ ab=-q $. Then,\\
$$ x^2-px-q=0 .$$
Since, both $\mbox{a} \& \mbox{b}$ are roots, they must satisfy the equation. Thus, the above equation becomes;\\
$$ a^2-pa-q=0 \qquad \mbox{and} \qquad b^2-pb-q=0. $$\\
Starting from $a^2-pa-q=0$, we get recursion for powers;\\

$a^2=pa+q(1)$ \\

$a^3=(p^2+q)a+qp$ \\

$a^4=p(p^2+2q)a+q(p^2+q)$ \\

$a^5=[p^2(p^2+2q)+q(p^2+q)]a+qp(p^2+2q).$ \\

Introducing sequence of polynomials of two variables $F_n(p,q)$; \\
$F_0(p,q)=0$\\
$F_1(p,q)=1$\\
$F_2(p,q)=p$\\
$F_3(p,q)=p^2+q$\\
$F_4(p,q)=p(p^2+2q)$\\
$F_5(p,q)=p^2(p^2+2q)+q(p^2+q)$\\
$F_6(p,q)=p^3(p^2+2q)+2qp(p^2+q)+q^2p$ \\
$F_7(p,q)=p^4(p^2+2q)+3p^4q+6p^2q^2+q^3$\\
\phantom{abc}.\\
\phantom{abc}.\\
we get $n^{th}$ power by recursion and find ;
\begin{equation}
a^n=F_n(p,q)a+qF_{n-1}(p,q) \label{nthpower}.\\
\end{equation}
By applying the same procedure for b, $b^n=F_n(p,q)b+qF_{n-1}(p,q)$ is obtained. Both formulas can be proved by induction.\\
By subtracting $b_n$ from $a_n$, we get Binet formula for these polynomials:\\
\begin{equation}
F_n{(p,q)}=\frac{a^n-b^n}{a-b}, \label{generalizedbinet}
\end{equation}\\
where $a,b=\frac{p}{2}\pm\sqrt{\frac{p^2}{4}+q}.$\\
This way, we get Binet type formula for sequence of Fibonacci polynomials.\\
Now, we are going to derive recursion formula for  Fibonacci polynomials;\\
$$F_{n+1}(p,q)=AF_{n}(p,q)+BF_{n-1}(p,q)$$\\
To find the coefficients A and B in above equation we substitute the Binet formulas for Fibonacci polynomials (\ref{generalizedbinet}) in $F_{n+1},F_{n},F_{n-1}$. Thus we have;\\
$$a^{n+1}=Aa^n+Ba^{n-1}$$\\
$$b^{n+1}=Ab^n+Bb^{n-1}$$
Equivalently;\\
$$a^2=Aa+B$$\\
$$b^2=Ab+B$$\\
Since we found that $a^2=pa+q$ \,(when we looked for the recursion relation for powers a), by substituting it;\\
$$a^2=Aa+B \; \longrightarrow pa+q=Aa+B \; \longrightarrow \; A=p \; \mbox{and} \; B=q.$$\\
Also, similarly from $b^2=Ab+B$ we can get the same result. Now we have,\\
\begin{equation}
F_{n+1}(p,q)=pF_{n}(p,q)+qF_{n-1}(p,q).
\end{equation}\\
Thus, we obtained recursion relation for sequence of Fibonacci polynomials. If p and q are arbitrary integer numbers, then we get the sequence of integer numbers;\\
$F_0(p,q)=0$\\
$F_1(p,q)=1$\\
$F_2(p,q)=p$\\
$F_3(p,q)=p^2+q$\\
$F_4(p,q)=p(p^2+q)+qp$\\
$F_5(p,q)=p^2(p^2+2q)+q(p^2+q)$\\
$F_6(p,q)=p^3(p^2+2q)+2qp(p^2+q)+q^2p$ \\
$F_7(p,q)=p^4(p^2+2q)+3p^4q+6p^2q^2+q^3$\\
\phantom{abc}.\\
\phantom{abc}.\\
which we call Fibonacci polynomial numbers.\\

Also;  when we choose $p=q=1$, the recursion relation will be standard recursion and Fibonacci numbers come. So; $F_n(1,1)=F_n.$

\subsection{Generalized Fibonacci Polynomials}

Fibonacci polynomials with initial values $F_0=0$\;\;and\;\;$F_1=1$ can be generalized to arbitrary initial values $G_0$  and $G_1$. So we define generalized Fibonacci polynomials $G_n(p,q)$ by the recursion formula
\begin{equation}
G_{n+1}(p,q)=pG_{n}(p,q)+qG_{n-1}(p,q), \label{gfp1}
\end{equation}
with initial values.
\begin{equation}
G_{0}(p,q)=G_{0},G_{1}(p,q)=G_{1}. \label{gfp2}
\end{equation}
It is easy to show that generalized Fibonacci polynomials can be represented as superposition of Fibonacci polynomial sequences:\\
\begin{equation}
G_{n}(p,q)=G_{1}F_{n}(p,q)+qG_{0}F_{n-1}(p,q). \label{gfp3}
\end{equation}
For generalized Fibonacci polynomials we find the following Binet type formula;\\

$$G_{n}(p,q)=\frac{(G_1-bG_0)a^n-(G_1-aG_0)b^n}{a-b},$$\\
where\;\; $a,b=\frac{p}{2}\pm\sqrt{\frac{p^2}{4}+q}.$\\

First few Generalized Fibonacci numbers are;\\
$G_0(p,q)=G_0$ \\
$G_1(p,q)=G_1$ \\
$G_2(p,q)=G_1p+qG_0$ \\
$G_3(p,q)=G_1(p^2+q)+qpG_0$ \\
$G_4(p,q)=G_1(p(p^2+q)+qp)+q(p^2+q)G_0$ \\
$G_5(p,q)=G_1(p^2(p^2+2q)+q(p^2+q))+q(p(p^2+q)+qp)G_0$ \\
\phantom{abc}.\\
\phantom{abc}.\\

If initial values $G_0=k$  and $G_1=l$ are integer, as well as $p=s$ and $q=t$ coefficients
then, we get Generalized Fibonacci polynomial numbers. \\

Also as a special case, if we choose $G_0=0$ and $G_1=p=q=1$, then we get the sequence of Fibonacci numbers. \\
\section{Applications}

\subsection{Fibonacci subsequences}

\qquad
If we consider a subsequence of the Fibonacci sequence, then this subsequence of numbers should satisfy some rules, which in general are different from the Fibonacci addition formula. Here we consider special type of subsequences generated by equidistant numbers, or the congruent numbers.

Let us consider the sequence $G_n$ as subsequence of Fibonacci numbers $G_n=F_{3n}$, which corresponds to equidistant integers $3,3+3,... ,3n$, or $n=0 \phantom{ab} (mod \phantom{.} 3)$. We like to know the recursion relation and corresponding initial conditions for this sequence.
For this, we have to find A and B coefficients in equation:
\begin{equation}
G_{n+1}=AG_{n}+BG_{n-1} \label{triplecoefcnt.}
\end{equation}
or equivalently,\quad $F_{3n+3}=AF_{3n}+BF_{3n-3}$.
By using Fibonacci recursion we rewrite $F_{3n+3}$ in terms of bases $F_{3n}$ and $F_{3n-3}$ as follows;

\begin{align}
\quad F_{3n+3}&=(F_{3n+2})+F_{3n+1} \nonumber \\
&=(F_{3n+1}+F_{3n})+F_{3n+1} \nonumber \\
&=2(F_{3n+1})+F_{3n} \nonumber\\
&=2(F_{3n}+F_{3n-1})+F_{3n} \nonumber\\
&=3F_{3n}+2(F_{3n-1}) \nonumber \\
&=3F_{3n}+2(F_{3n-2}+F_{3n-3})\nonumber \\
&=3F_{3n}+2F_{3n-3}+(2F_{3n-2})\nonumber
\end{align}

\begin{align}
\qquad \;\;\;\;\;\;\;\;\; &=3F_{3n}+2F_{3n-3}+(F_{3n}-F_{3n-3}) \nonumber \\
&=4F_{3n}+F_{3n-3}. \nonumber
\end{align}
As a result we find the recursion formula;
$$G_{n+1}=4G_{n}+G_{n-1}.$$
This recursion formula shows that this sequence is the Generalized Fibonacci polynomial sequence with initial values
$G_0=0$, $G_1=2$ and $p=4$, $q=1$.

\subsubsection{Equi-Fibonacci Sequences}
Here we introduce equidistant numbers determined by formula $kn+\alpha$, where $\alpha=0,1,2,...,k-1$ and $k=1,2,...$ are fixed numbers, and $n=0,1,2,3,...$. Distance between such numbers is k, this way we call them equi-distant numbers. The equidistant numbers are congruent numbers $\alpha$ (mod \phantom{.}k). Fibonacci numbers $F_{kn+\alpha}$ corresponding to such equi-distant numbers we call Equi-Fibonacci numbers. Now we are going to show that equi-Fibonacci number subsequence of Fibonacci sequence is Generalized Fibonacci Polynomial number sequence.\\
Let us consider subsequence of Fibonacci numbers determined by formula ${G^{(k;\alpha)}_n} \equiv F_{kn+\alpha}$. This subsequence satisfies the recursion formula according to the next theorem.\\

\textbf{THEOREM 1:}\quad Subsequence ${G^{(k;\alpha)}_n} \equiv F_{kn+\alpha}$ is subject to the next recursion formula
\begin{equation}
G^{(k;\alpha)}_{n+1}=L_k{G^{(k;\alpha)}_n}+(-1)^{k-1}{G^{(k;\alpha)}_{n-1}}, \label{section:Triple relations1}
\end{equation}
where $L_k$ are Lucas numbers. This is why it is given by Generalized Fibonacci polynomials with integer arguments as
\begin{equation}
G^{(k;\alpha)}_{n}=F_{kn+\alpha}=G_n{(L_k,(-1)^{k-1})}.
\end{equation}
\\
For proof see Appendix 7.1\\

As an example we apply formula (\ref{section:Triple relations1}), to the sequence $G^{(3;0)}_{n}=F_{3n}$,
$$G^{(3;0)}_{n+1}=L_3G^{(3;0)}_{n}+(-1)^{3-1}G^{(3;0)}_{n-1}=4G^{(3;0)}_{n}+G^{(3;0)}_{n-1}$$

and  find the same recursion relation as before,\\
$$G^{(3;0)}_{n+1}=4G^{(3;0)}_{n}+G^{(3;0)}_{n-1}.$$\\
Also, the same recursion is valid for both sequences $G^{(3;1)}_n=F_{3n+1}$, and $G^{(3;2)}_n=F_{3n+2}$.\\
As we can see from these sequences;\\
$G^{(3;0)}_n=F_{3n}=F_{0\phantom{..}( \mbox{mod}\phantom{..} 3)}=0,2,8,34,..$\\
$G^{(3;1)}_n=F_{3n+1}=F_{1\phantom{..}( \mbox{mod}\phantom{..} 3)}=1,3,13,55,..$\\
$G^{(3;2)}_n=F_{3n+2}=F_{2\phantom{..}( \mbox{mod}\phantom{..} 3)}=1,5,21,89,..$\\
each of three sequences starts with different initial values and they cover the whole Fibonacci sequence.\\

We can easily say that when we generate sequence of $F_{kn}$, the same recursion will be valid to generate the sequences $F_{kn+1}$,...,$F_{kn+(k-1)}$ and all sequences begin with different initial conditions. At the end they cover the whole Fibonacci sequence. \\

As an example, for even numbers $n=2m$  the subsequence of Fibonacci numbers $F_{2m}=G_n$ satisfies recursion formula $$G_{n+1}=3 G_{n}-G_{n-1}.$$\\

For odd numbers $n=2m+1$, the subsequence of Fibonacci numbers $F_{2m+1}=G_n$ satisfies the same recurrence relation.\\

For the case $k=3$, we have three subsequences of natural numbers as $n=3m, 3m+1, 3m+2$\,, then corresponding subsequences of Fibonacci numbers $F_{3m},\,F_{3m+1},\,F_{3m+2}$ satisfy the recursion formula $$G_{n+1}=4 G_{n}+G_{n-1}.$$\\

Case k; for the subsets of natural numbers $km,km+1,...,km+k-1$ we have subsequences of Fibonacci numbers satisfying recursion formula $$G_{n+1}=L_k G_{n}+(-1)^{k-1}G_{n+1}.$$\\
Below schema shows valid recursion relations for the desired sequences which have the $\alpha$ values, respectively;\\

\quad SEQUENCES \hspace{1cm} DIFFERENCE \hspace{1cm} VALID RECURSION RELATION \\
$$\quad \quad \hspace{0,5cm}   G^{(1,0)}_n=\{F_{n}\}  \quad \quad\hspace{1,5cm}  0  \quad \quad \hspace{1cm} G_{n+1}^{(1;0)}=G_{n}^{(1;0)}+G_{n-1}^{(1;0)}\hspace{0,25cm}  \;k=1,\alpha=0$$ \\
$$\qquad  G^{(2,\alpha)}_n=\{F_{2n},F_{2n+1}\}  \quad \hspace{1,05cm}1 \quad \quad \quad \quad G_{n+1}^{(2;\alpha)}=3G_{n}^{(2;\alpha)}-G_{n-1}^{(2;\alpha)}\hspace{0,5cm}  \;k=2,\alpha=0,1$$ \\
$$\,\,G^{(3,\alpha)}_n=\{{F_{3n},F_{3n+1},F_{3n+2}}\} \hspace{1,05cm}2 \quad \quad \quad G_{n+1}^{(3;\alpha)}=4G_{n}^{(3;\alpha)}+G_{n-1}^{(3;\alpha)}\hspace{0,5cm}  \;k=3,\alpha=0,1,2$$ \\
\phantom{abc}.\\
\phantom{abc}.\\$$G^{(k,\alpha)}_n=\{{F_{kn},...,F_{kn+(k-1)}}\}\quad \hspace{0,5cm}k-1\qquad G_{n+1}^{(k;\alpha)}=L_kG_{n}^{(k;\alpha)}+(-1)^{k-1}G_{n-1}^{(k;\alpha)}\hspace{0,3cm}\;k,\alpha=0...k-1$$\\

\subsubsection{Equi-Fibonacci superposition}
Equi-Fibonacci numbers are determined in terms of Generalized Fibonacci polynomials as follows,
\begin{equation}
G^{(k;\alpha)}_{n}=F_{kn+\alpha}=G_n{(L_k,(-1)^{k+1})} \label{super}
\end{equation}
$$G_{0}^{(k,\alpha)}=F_\alpha \;\;,\;\; G_{1}^{(k,\alpha)}=F_{\alpha+k}$$\\
Then, by using superposition formula for Generalized Fibonacci Polynomials (\ref{gfp3}) we get equi-Fibonacci numbers as superposition
$$G^{(k;\alpha)}_{n}=F_{kn+\alpha}=F_{\alpha+k}F_n^{(k)}+F_{\alpha}(-1)^{k+1}F_{n-1}^{(k)}$$ of Higher Order Fibonacci numbers. Easy calculation shows that Higher Order Fibonacci numbers can be written as ratio;\\
$$F_n^{(k)}=\frac{(\varphi^k)^n-(\varphi'^k)^n}{\varphi^k-\varphi'^k}=\frac{F_{nk}}{F_k}.$$
This Higher Order Fibonacci numbers satisfy the same recursion formula as (\ref{section:Triple relations1}).
\subsection{Limit of Ratio Sequences}

Binet Formula for Fibonacci numbers allows one to find the limit of ratio sequence $U_n=\frac{F_{n+1}}{F_{n}}$ as a Golden Ratio. Here, we will derive similar limits for Generalized Fibonacci numbers and Fibonacci Polynomials.

\subsubsection{Fibonacci Sequence and Golden Ratio}

\quad Golden Ratio appears many times in geometry, art, architecture and other areas to analyze the proportions of natural objects. \\
There is a hidden beauty in Fibonacci sequence. While going to infinity, if we take any two successive Fibonacci numbers, their ratio become more and more close to the Golden Ratio $"\varphi"$, which is approximately 1,618.\\
In a mathematical way we can see it as;\\
$$\mbox{Since} \quad \varphi'=-\frac{1}{\varphi} , \;\lim_{n{\to \infty}} \varphi'^{n+1}=\lim_{n{\to \infty}} \varphi'^{n}=0,$$\\
then,
$$\lim_{n{\to \infty}} \frac{F_{n+1}}{F_{n}}=\lim_{n{\to \infty}} \frac{\varphi^{n+1}-\varphi'^{n+1}}{\varphi^{n}-\varphi'^{n}}=\lim_{n{\to \infty}} \frac{\varphi^{n+1}-\frac{1}{\varphi^{n+1}}}{\varphi^{n}-\frac{1}{\varphi^{n}}}     =\lim_{n{\to \infty}} \frac{\varphi^{n+1}}{\varphi^{n}}=\varphi.$$ \\

\subsubsection{Two Fibonacci subsequences}
$$G_n=\frac{(G_1-\varphi'G_0)\varphi^n-(G_1-\varphi G_0)\varphi'^n}{\varphi-\varphi'},$$
$$\lim_{n{\to \infty}} \frac{G_{n+1}}{G_n}=\varphi.$$

\subsubsection{Arbitrary number of subsequences}
$$G_n^{(k)}=\frac{\varphi^nP_k(\varphi)-\varphi'^nP_k(\varphi')}{\varphi-\varphi'},$$
$$\lim_{n{\to \infty}} \frac{G^{(k)}_{n+1}}{G^{(k)}_n}=\varphi.$$

\subsubsection{Fibonacci polynomials}
From Binet formula  for Fibonacci polynomials;\\
$$F_{n}{(a,b)}=\frac{a^n-b^n}{a-b},\,\quad\mbox{where}\,\quad a,b=\frac{p}{2}\pm\sqrt{\frac{p^2}{4}+q}$$\\
we get the limit,\\
$$\lim_{n{\to \infty}} \frac{F_{n+1}{(a,b)}}{F_{n}{(a,b)}}=\mbox{max}(a,b)=a.$$

\section{Binet-Fibonacci Curve}

\subsection{Fibonacci Numbers of Real argument}
\quad In Binet formula if we consider index t as an integer, corresponding Fibonacci numbers will be integer numbers,\\
 $$ t \in \mathbb{Z}  \rightarrow F_{t} \in \mathbb{Z}. $$ \\
When we choose t as a real number, corresponding $F_t$ will be in complex plane \cite{PN}. Thus, we have;\\
 $$t \in \mathbb{R}  \rightarrow  F_{t} \in \mathbb{C}. $$ \\
Let us analyze the second case;\\
$$F_t=\frac{\varphi^t-\varphi'^t}{\varphi-\varphi'}, \quad \mbox{where} \; \mbox{t} \in \mathbb{R}.$$\\
In this formula, we have now multi-valued function;
$$(-1)^t=e^{t log(-1)}=e^{i \pi t(2n+1)},$$ where $n=0,\pm1,\pm2,...$
In following calculations we choose only one branch of this function with $n=0:$\\
$$\varphi'^t=\left(\frac{-1}{\varphi}\right)^t=\frac{(-1)^t}{\varphi^t}=\frac{e^{i \pi t}}{\varphi^t},$$\\
\begin{align}
F_{t}=&\frac{\varphi^t-(\frac{-1}{\varphi})^t}{\varphi-\varphi'} \nonumber \\
F_{t}=&\frac{e^{t\ln(\varphi)}-e^{i \pi t}\varphi^{-t}}{\varphi-(\frac{-1}{\varphi})} \nonumber \\
F_{t}=&\frac{1}{\varphi+\varphi^{-1}}\,\,\,\left(e^{t\ln(\varphi)}-e^{i \pi t}\,\,\,e^{-t\ln(\varphi)}\right) \nonumber \\
F_{t}=&\frac{1}{\varphi+\varphi^{-1}}\,\,\,\left(e^{t\ln(\varphi)}-(\cos(\pi t)+i.\sin(\pi t))e^{-t\ln(\varphi)}\right) \nonumber \\
F_{t}=&\frac{1}{\varphi+\varphi^{-1}}\,\,\,\left(e^{t\ln(\varphi)}-\cos(\pi t)e^{-t\ln(\varphi)}-i\sin(\pi t))e^{-t\ln(\varphi)}\right) \nonumber \\
F_{t}=&\frac{1}{\varphi+\varphi^{-1}}\,\,\,(e^{t\ln(\varphi)}-\cos(\pi t)e^{-t\ln(\varphi)})+i\frac{1}{\varphi+\varphi^{-1}}(-\sin(\pi t))e^{-t\ln(\varphi)}) \nonumber
\end{align} \\

This way, we have complex valued function $F_{t}$ of real variable t. This function describes a curve in complex plane parameterized by real variable t. This curve for -$\infty$$<$t$<$$\infty$, we call the ``Binet-Fibonacci curve". Parametric form of this curve is;\\

$$F_{t}=\left(x(t),y(t)\right)\,;$$\\
$$Re(F_t)=x(t)=\frac{1}{\varphi+\varphi^{-1}}\left(e^{t\ln(\varphi)}-\cos(\pi t)e^{-t\ln(\varphi)}\right)$$\\
$$Im(F_t)=y(t)=\frac{1}{\varphi+\varphi^{-1}}\left(-\sin(\pi t)e^{-t\ln(\varphi)}\right)$$\\

We plot this curve for 0$<$t$<$$\infty$, in Figure 1 by Mathematica.

\begin{center}\vspace{0.25cm}
\includegraphics[width=1\linewidth]{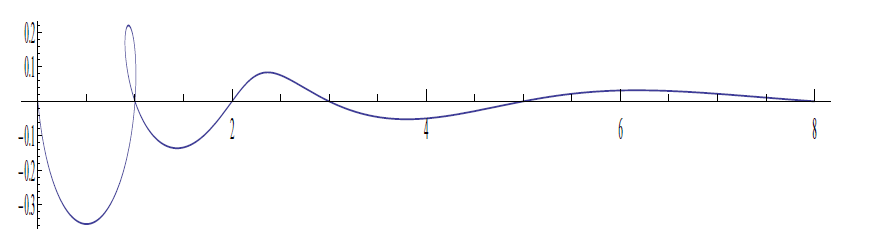}
\captionof{figure}{Binet-Fibonacci Oscillating Curve (B.F.O.C.)}
\end{center}\vspace{0.25cm}

The curve intersects real axis x at Fibonacci numbers $F_n$ corresponding to integer values of parameter t=n. \\
Below we discuss some properties of this curve.

\subsection{Area Sequences for Binet-Fibonacci\\ Oscillating Curve}

\quad Here we analyze the area sequences under Binet-Fibonacci Oscillating curve. We use the Green's Theorem area formula to calculate the area value of segments;\\

$$A_{n,n+1}=\frac{1}{2} \int_n^{n+1} (x dy-y dx).$$ \\
By rewriting $dy=\left(\frac{dy}{dt}\right)dt$ and $dx=\left(\frac{dx}{dt}\right)dt$, our formula becomes;\\
\begin{equation}
A_{n,n+1}=\frac{1}{2}\int_n^{n+1} (x \dot{y}-y \dot{x})dt. \label{greensthmareaformula}\\
\end{equation}
Since we know $x(t)$ \& $y(t)$ as components of Binet formula for real argument, we substitute them into equation (\ref{greensthmareaformula}) and by taking integral we obtain the area formula for the desired segments of B.F.O.C. We put the - sign in front of the area formula to be hold with signs of calculated area segments;\\
\begin{equation}
A_{n,n+1}=-\frac{1}{10} \left[\frac{4(-1)^nln(\varphi)}{\pi}-\frac{\pi}{2ln(\varphi)}(F_{2n}-\frac{1}{\varphi}F_{2n+1})\right]. \\
\end{equation}

From this formula, we have interesting result in B.F.O.C, when we look the limit value of area segment at infinity;\\

\begin{align} \nonumber
\lim_{n{\to \infty}} A_{n,n+1}&=\lim_{n{\to \infty}}  -\frac{1}{10}\left[\frac{4(-1)^nln(\varphi)}{\pi}-\frac{\pi}{2ln(\varphi)} \varphi'^{2n+1})\right] \nonumber \\
&=\lim_{n{\to \infty}} -\frac{1}{10}\left[\frac{4(-1)^nln(\varphi)}{\pi}-\frac{\pi}{2ln(\varphi)}(-\frac{1}{\varphi})^{2n+1})\right] \nonumber \\
&=\lim_{n{\to \infty}} -\frac{1}{10}\left[\frac{4(-1)^nln(\varphi)}{\pi}+\frac{\pi}{2ln(\varphi)}\frac{1}{\varphi^{2n+1}}\right] \nonumber \\
 \nonumber
\end{align}
Thus, when n goes to infinity we find that the sequence of segment's area has finite limit;\\
\begin{equation}
A_\infty=\lim_{n{\to \infty}} |A_{n,n+1}|=\left|-\frac{1}{10}\frac{4(-1)^nln(\varphi)}{\pi}\right|=\frac{2}{5\pi}\,ln(\varphi). \label{area}
\end{equation}\\

Since we have \;$A_\infty=\frac{2}{5\pi}\,ln(\varphi)$, other possible identities can also be written; \\

$$\pi A_\infty=\frac{2}{5}\,ln(\varphi)\quad \mbox{and} \quad \varphi=e^{\frac{5\pi}{2}A_\infty}$$\\

We can say that at infinity, the area of segments becomes close to the value $\frac{2}{5}\frac{ln(\varphi)}{\pi}\approx 0.06\;.$\\

Remarkable is that result (\ref{area}) includes three different fundamental constants of mathematics as $\varphi,\pi$ and $e$.

As we can see from B.F.O.C, area segments starting with even number have negative sign, and area segments starting with odd number have positive sign as an area value. According to this observation we formulate our next theorem.\\

\textbf{THEOREM 2:}\quad$A_{n,n+1}<0$ if $n=2k$ and $A_{n,n+1}>0$ if $n=2k+1$, where k=1,2,3,...\\

For proof see Appendix 7.2\;.\\

If we look to the sign of addition two consecutive area segments under B.F.O.C, we obtain interesting observation.\\

\textbf{THEOREM 3:}\quad $A_{n,n+1}+A_{n+1,n+2}<0$ for both $n=2k$ and $n=2k+1$ n=2,3,4,..\\

For proof see Appendix 7.2\;.

 In Binet formula for real argument, segment's area formula was;\\

 $$A_{n,n+1}=-\frac{1}{10} \left[\frac{4(-1)^nln(\varphi)}{\pi}-\frac{\pi}{2ln(\varphi)}(F_{2n}-\frac{1}{\varphi}F_{2n+1})\right].$$\\

 Also we can write;\\

$$A_{n-1,n}=-\frac{1}{10} \left[\frac{4(-1)^{n-1}ln(\varphi)}{\pi}-\frac{\pi}{2ln(\varphi)}(F_{2n-2}-\frac{1}{\varphi}F_{2n-1})\right].$$\\

Here, we define $\lambda_n$ as ratio of areas;\\

$$\lambda_n=\left|\frac{A_{n,n+1}}{A_{n-1,n}}\right|=\left|\frac{8(ln(\varphi))^2(-1)^n-\pi^2(F_{2n}-\frac{1}{\varphi}F_{2n+1})}{-8(ln(\varphi))^2(-1)^n-\pi^2(F_{2n-2}-\frac{1}{\varphi}F_{2n-1})}\right|$$\\

and after some calculations (see Appendix 7.3), as a result;\\
$$\;\lambda_n=\left|\frac{A_{n,n+1}}{A_{n-1,n}}\right|=\frac{1}{\varphi^2}\frac{|\frac{8(ln(\varphi))^2}{\pi^2}\varphi^{2n}\varphi+(-1)^n|}{|\frac{8(ln(\varphi))^2}{\pi^2}\varphi^{2n}\frac{1}{\varphi}-(-1)^n|}\;\;\;\;\;\;\mbox{is obtained.}$$

If $n=2k$;\\

$$\lambda_{2k}=\left|\frac{A_{2k,2k+1}}{A_{2k-1,2k}}\right|=\frac{1}{\varphi^2}\frac{|\frac{8(ln(\varphi))^2}{\pi^2}\varphi^{4k+1}+1|}{|\frac{8(ln(\varphi))^2}{\pi^2}\varphi^{4k-1}-1|}=\frac{|\frac{8(ln(\varphi))^2}{\pi^2}\varphi^{4k-1}+\frac{1}{\varphi^2}|}{|\frac{8(ln(\varphi))^2}{\pi^2}\varphi^{4k-1}-1|}>1.$$\\

If $n=2k+1$;\\

$$\lambda_{2k+1}=\left|\frac{A_{2k+1,2k+2}}{A_{2k,2k+1}}\right|=\frac{1}{\varphi^2}\frac{|\frac{8(ln(\varphi))^2}{\pi^2}\varphi^{4k+3}-1|}{|\frac{8(ln(\varphi))^2}{\pi^2}\varphi^{4k+1}+1|}=\frac{|\frac{8(ln(\varphi))^2}{\pi^2}\varphi^{4k+1}-\frac{1}{\varphi^2}|}{|\frac{8(ln(\varphi))^2}{\pi^2}\varphi^{4k+1}+1|}<1.$$\\
Thus Theorem 2, which is about sign of the addition two consecutive segments, is proved in other logical way.\\

Now, let us look;\\

$$\lim_{n{\to \infty}} \lambda_n=\frac{1}{\varphi^2}\frac{|\varphi|}{|\frac{1}{\varphi}|}=\frac{1}{\varphi^2}\varphi^2=1.$$\\

From this calculation, we can interpret that at infinity, the area values of consecutive segments are equal. In fact, we found before that at infinity, segment's area value is $\frac{2}{5\pi}\,ln(\varphi)$ under Binet-Fibonacci Oscillating curve. In a logical manner, the ratio of consecutive area segments will be equal at infinity.\\

Now, we calculate the curvature of the Binet-Fibonacci Oscillating curve. It is calculated by using the formula;\\
$$\kappa(t)=\frac{|\dot{r}\times\ddot{r}|}{|\dot{r}|^3}$$\\
where,
$$\vec{r}(t)=F_{t}=(x(t),y(t))=\frac{1}{\varphi+\varphi^{-1}}\left(e^{t\ln(\varphi)}-\cos(\pi t)e^{-t\ln(\varphi)},-\sin(\pi t)e^{-t\ln(\varphi)}\right)$$\\

Then, by substituting vector $\vec{r}$, we find how much the Binet-Fibonacci Oscillating curve is curved at any given point t;\\

\begin{equation}
\kappa(t)=\frac{3\pi ln^{2}\varphi \cos(\pi t)+\sin(\pi t)(\pi^{2}ln\varphi-2ln^{3}\varphi)+e^{-2tln\varphi}(\pi(\pi^2+ln^2\varphi))}{5\left[\frac{1}{5}\left(ln^{2}\varphi(e^{2tln\varphi}+e^{-2tln\varphi})+2ln\varphi\; \pi\sin(\pi t)+2ln^{2}\varphi\;\cos(\pi t)+\pi^{2}e^{-2tln\varphi}\right)\right]^{3/2}} \nonumber \\
\end{equation}\\

By using this result, we find that;\\
\begin{equation}
\lim_{n{\to \infty}} \kappa(n)=0. \nonumber \\
\end{equation}
It means, when n is close to infinity, at Fibonacci number points the curvature is 0. Namely, the curve behaves like line at Fibonacci numbers.\\

Another amazing result is found when we look the curvature ratio at consecutive Fibonacci numbers;\\

$$\lim_{n{\to \infty}} \frac{\kappa(n+1)}{\kappa(n)}=-\frac{1}{\varphi^3}.$$\\
In this result,  ``-" sign comes from its signed curvature of the curve.\\

By using this result we can also calculate;\\
$$\lim_{n{\to \infty}} \frac{\kappa(n+2)}{\kappa(n)}=\lim_{n{\to \infty}} \frac{\kappa(n+2)}{\kappa(n+1)}. \frac{\kappa(n+1)}{\kappa(n)}=\left(-\frac{1}{\varphi^3}\right)\left(-\frac{1}{\varphi^3}\right)=\frac{1}{\varphi^6}.$$

We can generalize this result as;
$$\lim_{n{\to \infty}} \frac{\kappa(n+k)}{\kappa(n)}=\left(-\frac{1}{\varphi^3}\right)^k.$$
\newpage
\subsection{Fibonacci Spirals}
\quad Up to now, we have thought only positive real numbers $t$ as an argument in Binet formula for $F_t$. But now, we look what will be if we consider t as negative real numbers: -$\infty$$<$t$<$0\\
In equation (\ref{relationwith negative index}), we have found the relation between $F_n$ and $F_{-n}$ so that \\
$$F_{-n}=(-1)^{n+1}F_{n}.$$\\
\phantom{ab}From this equality, we can see that for negative even number n, the values of the function are given by negative Fibonacci numbers. For odd numbers n we will have positive value for Fibonacci numbers with odd index.\\
\phantom{ab}In contrast to positive t when function $F_t$ is oscillating along positive half axis, where $Re(F_t)$ is positive, for negative t we are getting a spiral, crossing real line from both sides. We called such curve as ``Binet-Fibonacci Spiral curve" since it intersects real axis at Fibonacci numbers. Its shape in complex plane is obtained by using parametric plot in Mathematica, in Figure 2.
\begin{center}\vspace{0.25cm}
\includegraphics[width=0.7\linewidth]{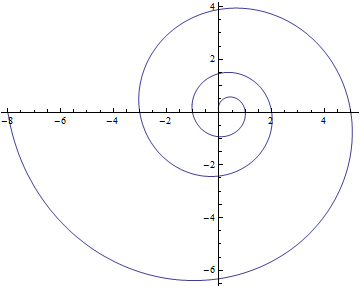}
\captionof{figure}{Binet-Fibonacci Spiral Curve}
\end{center}\vspace{0.25cm}


\begin{center}\vspace{0.25cm}
\includegraphics[width=0.90\linewidth]{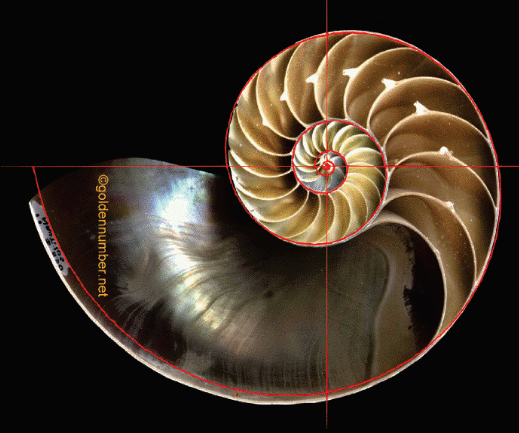}
\captionof{figure}{Nautilus 1}
\end{center}\vspace{0.25cm}


\begin{center}\vspace{0.25cm}
\includegraphics[width=0.40\linewidth]{spiral}
\hspace{1.5cm}
\includegraphics[width=0.40\linewidth]{nautilus1}
\captionof{figure}{Comparison of Binet-Fibonacci Spiral Curve with Nautilus 1}
\end{center}\vspace{0.25cm}


\begin{center}\vspace{0.25cm}
\includegraphics[width=0.90\linewidth]{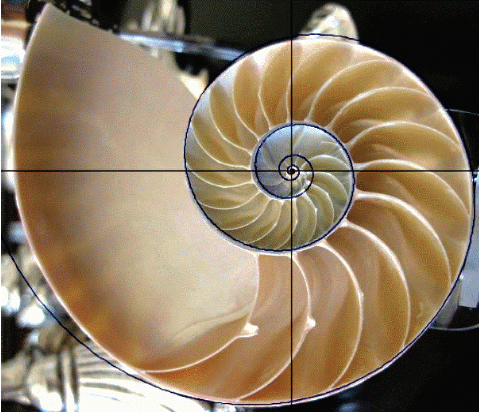}
\captionof{figure}{Nautilus 2}
\end{center}\vspace{0.25cm}


\begin{center}\vspace{0.25cm}
\includegraphics[width=0.40\linewidth]{spiral}
\hspace{1.5cm}
\includegraphics[width=0.40\linewidth]{nautilus2}
\captionof{figure}{Comparison of Binet-Fibonacci Spiral Curve with Nautilus 2}
\end{center}\vspace{0.25cm}

To see that there is no other intersection points with the x-axis, apart from Fibonacci numbers, let us look $Im(F_t).$\\

When $Im(F_t)=0$, we get intersection points with x-axis as,\\
$$Im(F_t)=\frac{1}{\varphi+\varphi^{-1}}\left(-\sin(\pi t)e^{-t\ln(\varphi)}\right)=0.$$\\
Since $e^{-t\ln(\varphi)}\neq0$,\\
$$\sin(\pi t)=0 \Rightarrow \pi t=\pi k \Rightarrow\;\; t=k, \;\;\mbox{where}\;\; k=0,\pm1,\pm2,\pm3,... $$\\
\phantom{abc}This means at points $F_k$, where $k=0,\pm1,\pm2,\pm3,...$\,\,, we have intersection points with real axis. So, we can say that there is no other intersection points of the Fibonacci spiral with the real axis, except at Fibonacci numbers $F_k$.\\

Using the result, on curvature of Binet-Fibonacci Oscillating curve, we can expect that when $t\rightarrow-\infty$, $\kappa(t)\rightarrow0.$ Intuitively, due to the having huge radius, it is clear that Binet-Fibonacci Spiral behaves like line at any given point when $t\rightarrow-\infty$. \\

In addition, we look the curvature ratio when $t{\to -\infty}$;\\

$$\lim_{t{\to -\infty}} \frac{\kappa(t+1)}{\kappa(t)}=\frac{1}{\varphi}.$$\\

Furthermore, we can generalize this result;\\
$$\lim_{t{\to -\infty}} \frac{\kappa(t+k)}{\kappa(t)}=\lim_{t{\to -\infty}} \frac{\kappa(t+k)}{\kappa(t+(k-1))}\;.\;.\;.\;\frac{\kappa(t+1)}{\kappa(t)}=\left(\frac{1}{\varphi}\right)^k.$$\\

When we look at the Binet-Fibonacci spiral curve;\\
$$\lim_{n{\to -\infty}} \frac{F_{n+2}}{F_n}=\lim_{n{\to -\infty}} \frac{F_{n+2}}{F_{n+1}}\; \frac{F_{n+1}}{F_n}=\varphi^2=\varphi+1\approx 2.618\;.$$\\

Also, we approach to this value by looking along the x and y axis in our own Nautilus necklace (Figures 3 and 5).\\

When we compare our Fibonacci spiral with the Nautilus's natural spiral in nature, we realize that they have quite similar behaviours (Figures 4 and 6).\\

Also, we noticed that when $t\rightarrow-\infty$, our Binet-Fibonacci spiral curve behaves as a Logarithmic spiral. Furthermore, we know that at every point in logarithmic spiral we have constant angle between vector $\vec{r}$ and $\frac{d\vec{r}}{dt}$. We find that at infinity angle between the vector $\vec{r}$ and $\frac{d\vec{r}}{dt}$ is constant in the Binet-Fibonacci Spiral curve (asymptotically).\\

Finally; in Figure 7, we show the curve for both positive and negative real values of argument,-$\infty$$<$t$<$$\infty$, of the Binet Formula in the range of t between [-4,4].\\


\begin{center}\vspace{0.25cm}
\includegraphics[width=0.7\linewidth]{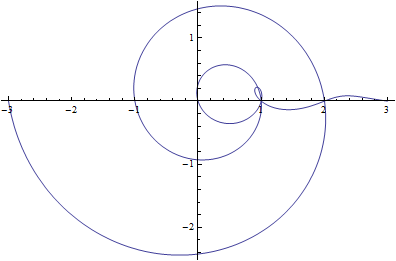}
\captionof{figure}{Binet-Fibonacci Curve}
\end{center}\vspace{0.25cm}


\section{Conclusions}

In the present paper, we have studied several generalizations of Fibonacci sequences. First, the generalized Fibonacci sequences as the sequences with arbitrary initial values. Second, the addition of two and more Fibonacci subsequences and Fibonacci polynomials with arbitrary bases. In all these cases we get Binet representation and asymptotics of the sequence as the ratio of consecutive terms. For Fibonacci numbers with congruent indices we derived general formula in terms of generalized Fibonacci polynomials and Lucas numbers. By extending Binet formula to arbitrary real variables we constructed Binet-Fibonacci curve in complex plane. For positive numbers the curve is oscillating with decreasing amplitude, while for negative numbers it becomes Binet-Fibonacci spiral. For t$<<$0, Binet-Fibonacci Spiral curve becomes Logarithmic spiral. Comparison with nautilus curve shows quite similar behavior. Areas under the curve and curvature characteristics are calculated as well as the asymptotic of relative characteristics.

\section{Appendix}

\subsection{Fibonacci Subsequence Theorem}

\textbf{THEOREM 1:}\quad Subsequence ${G^{(k;\alpha)}_n} \equiv F_{kn+\alpha}$ is subject to the next recursion formula
\begin{equation}
G^{(k;\alpha)}_{n+1}=L_k{G^{(k;\alpha)}_n}+(-1)^{k-1}{G^{(k;\alpha)}_{n-1}},
\end{equation} where $L_k$ are Lucas numbers. This is why it is given by Generalized Fibonacci polynomials with integer arguments as
\begin{equation}
G^{(k;\alpha)}_{n}=F_{kn+\alpha}=G_n{(L_k,(-1)^{k-1})}.
\end{equation}
\\
\textbf{PROOF 1:}\quad\\
To prove the equality ${G^{(k;\alpha)}_{n+1}}=L_k{G^{(k;\alpha)}_n}+(-1)^{k-1}{G^{(k;\alpha)}_{n-1}}$, we begin with left hand side of equation and at the end, we get right hand side of this equation.\\
Since, we can write ${G^{(k;\alpha)}_{n+1}}=F_{k(n+1)+\alpha}=F_{kn+k+\alpha}$;\\
\begin{align}
F_{kn+k+\alpha}&=\frac{1}{\varphi-\varphi'}\left[\varphi^{kn+k+\alpha}-\varphi'^{kn+k+\alpha}\right] \nonumber \\
&=\frac{1}{\varphi-\varphi'}\left[\varphi^{kn+\alpha}\varphi^k-\varphi'^{kn+\alpha}\varphi'^k \right] \nonumber \\
&=\frac{1}{\varphi-\varphi'}\left[\varphi^{kn+\alpha}\varphi^k+(-\varphi'^{kn+\alpha}\varphi^k+\varphi'^{kn+\alpha}\varphi^k)-\varphi'^{kn+\alpha}\varphi'^k\right]\nonumber\\
&=\frac{1}{\varphi-\varphi'}\left[(\varphi^{kn+\alpha}-\varphi'^{kn+\alpha})\varphi^k+\varphi'^{kn+\alpha}\varphi^k-\varphi'^{kn+\alpha}\varphi'^k\right]\nonumber\\
&=F_{kn+\alpha}\varphi^k+\frac{1}{\varphi-\varphi'}\left[\varphi'^{kn+\alpha}\varphi^k-\varphi'^{kn+\alpha}\varphi'^k\right] \nonumber\\
&=F_{kn+\alpha}(\varphi^k+(-\varphi'^k+\varphi'^k))+\frac{1}{\varphi-\varphi'}\left[\varphi'^{kn+\alpha}\varphi^k-\varphi'^{kn+\alpha}\varphi'^k\right]\nonumber\\
&=F_{kn+\alpha}(\varphi^k+\varphi'^k)-F_{kn+\alpha}\varphi'^k+\frac{1}{\varphi-\varphi'}\left[\varphi'^{kn+\alpha}\varphi^k-\varphi'^{kn+\alpha}\varphi'^k\right]\nonumber\\
&=F_{kn+\alpha}L_k+\frac{1}{\varphi-\varphi'}\left[\varphi'^{kn+\alpha}\varphi'^k-\varphi^{kn+\alpha}\varphi'^k+\varphi'^{kn+\alpha}\varphi^k-\varphi'^{kn+\alpha}\varphi'^k\right]\nonumber\\
&=F_{kn+\alpha}L_k+\frac{1}{\varphi-\varphi'}\left[\varphi'^{kn+\alpha}\varphi^k-\varphi^{kn+\alpha}\varphi'^k\right]\nonumber\\
&=F_{kn+\alpha}L_k+\frac{\varphi^k\varphi'^k}{\varphi-\varphi'}\left[\varphi'^{kn-k+\alpha}-\varphi^{kn-k+\alpha}\right]\nonumber\\
\nonumber
\end{align}
\begin{align}
&=F_{kn+\alpha}L_k-\frac{\varphi^k\varphi'^k}{\varphi-\varphi'}\left[\varphi^{kn-k+\alpha}-\varphi'^{kn-k+\alpha}\right]\nonumber\\
&=F_{kn+\alpha}L_k-\frac{(\varphi\varphi')^k}{\varphi-\varphi'}\left[\varphi^{kn-k+\alpha}-\varphi'^{kn-k+\alpha}\right] \;\; \mbox{since} \;\; (\varphi\varphi')^k=(-1)^k,\nonumber\\
&=F_{kn+\alpha}L_k-(-1)^k\left[\frac{\varphi^{kn-k+\alpha}-\varphi'^{kn-k+\alpha}}{{\varphi-\varphi'}}\right]\nonumber\\
&=F_{kn+\alpha}L_k+(-1)^{k+1}F_{kn-k+\alpha}\;\; \mbox{and since}\;\;  (-1)^{k+1}(-1)^{-2}=(-1)^{k-1}, \nonumber\\
&=F_{kn+\alpha}L_k+(-1)^{k-1}F_{kn-k+\alpha}\nonumber\\
\nonumber
\end{align}
Therefore, $F_{kn+k+\alpha}=F_{kn+\alpha}L_k+(-1)^{k-1}F_{kn-k+\alpha}$ is obtained. Equivalently, under our notation;\\
$${G^{(k;\alpha)}_{n+1}}=L_k{G^{(k;\alpha)}_n}+(-1)^{k-1}{G^{(k;\alpha)}_{n-1}}$$ is written.\\

\qquad \qquad \qquad \qquad  \qquad  \qquad \qquad  \qquad  \qquad \qquad  \qquad  \qquad  \qquad  $\blacksquare$\\

\subsection{Theorems on Alternating Sign of Areas }

\textbf{THEOREM 2:}\quad$A_{n,n+1}<0$ if $n=2k$ and $A_{n,n+1}>0$ if $n=2k+1$ where k=1,2,3,...\\

\textbf{PROOF 2:}\quad Firstly; if $n=2k$, then we have to show that;\\

$$A_{2k,2k+1}=-\frac{1}{10}\left[\frac{4ln(\varphi)}{\pi}-\frac{\pi}{2ln(\varphi)}(F_{4k}-\frac{1}{\varphi}F_{4k+1})\right]<0.$$

This inequality is valid if,\\

$$\frac{4ln(\varphi)}{\pi}>\frac{\pi}{2ln(\varphi)}(F_{4k}-\frac{1}{\varphi}F_{4k+1})$$\\

By induction we have; \quad $F_{4k}+\varphi'F_{4k+1}=\left(\varphi'\right)^{4k+1}.$ So, we have;\\

$$\frac{4ln(\varphi)}{\pi}>\frac{\pi}{2ln(\varphi)}(\varphi'^{4k+1})$$\\
$$\frac{8(ln(\varphi))^2}{\pi^2}>\left(\varphi'\right)^{4k+1}=\left(\frac{-1}{\varphi}\right)^{4k+1}=\left(\frac{1}{\varphi}\right)^{4k+1}(-1)^{4k+1}=-\left(\frac{1}{\varphi}\right)^{4k+1}$$\\So;
$$\frac{8(ln(\varphi))^2}{\pi^2}>-\left(\frac{1}{\varphi}\right)^{4k+1}$$\\
Thus, it is evident.\\

Secondly; if $n=2k+1$, we have to show that,\\
$$A_{2k+1,2k+2}=-\frac{1}{10} \left[-\frac{4ln(\varphi)}{\pi}-\frac{\pi}{2ln(\varphi)}(F_{4k+2}-\frac{1}{\varphi}F_{4k+3})\right]>0$$\\

It is valid when ;\\
$$-\frac{1}{10} \left[-\frac{4ln(\varphi)}{\pi}-\frac{\pi}{2ln(\varphi)}\left(\varphi'\right)^{4k+3}\right]>0$$\\

$$-\frac{1}{10} \left[-\frac{4ln(\varphi)}{\pi}+\frac{\pi}{2ln(\varphi)}\left(\frac{1}{\varphi}\right)^{4k+3}\right]>0$$\\

$$-\frac{1}{10}\left(\frac{\pi}{2ln(\varphi)}\right)\left[\frac{1}{(\varphi)^{4k+3}}-\frac{8(ln(\varphi))^2}{\pi^2}\right]>0$$\\
Now our next step is to prove;\\

$$\frac{8(ln(\varphi))^2}{\pi^2}>\frac{1}{\varphi^{4k+3}}$$\\

In fact, $$\frac{8(ln(\varphi))^2}{\pi^2}>\frac{1}{\varphi^7}>\frac{1}{\varphi^{4k+3}}\quad \mbox{for} \;\;k>1$$\\

So, we can easily say that $\frac{8(ln(\varphi))^2}{\pi^2}>\frac{1}{\varphi^{4k+3}}.$\\

As a result, proof is completed.\\

\qquad \qquad \qquad \qquad \qquad \qquad \qquad \qquad \qquad \qquad \qquad \qquad \qquad $\blacksquare$\\

\textbf{THEOREM 3:}\quad $A_{n,n+1}+A_{n+1,n+2}<0$ for both $n=2k$ and $n=2k+1$ n=2,3,4,..\\

\textbf{PROOF 3:}\quad Since we know; \\
$$A_{n,n+1}=-\frac{1}{10} \left[\frac{4(-1)^nln(\varphi)}{\pi}-\frac{\pi}{2ln(\varphi)}(F_{2n}-\frac{1}{\varphi}F_{2n+1})\right],$$\\
we can write addition of two consecutive area segments;\\

$A_{n,n+1}+A_{n+1,n+2}=$\\

$=-\frac{1}{10}\left[\frac{4(-1)^nln(\varphi)}{\pi}+\frac{4(-1)^{n+1}ln(\varphi)}{\pi}-\frac{\pi}{2ln(\varphi)}(F_{2n}-\frac{1}{\varphi}F_{2n+1})-\frac{\pi}{2ln(\varphi)}(F_{2n+2}-\frac{1}{\varphi}F_{2n+3})\right]$\\

$=\frac{1}{10}\left(\frac{\pi}{2ln(\varphi)}\right)\left[F_{2n}-\frac{1}{\varphi}F_{2n+1}+F_{2n+2}-\frac{1}{\varphi}F_{2n+3}\right]$\\

$=\frac{\pi}{20ln(\varphi)}\left[F_{2n}+F_{2n+2}-\frac{1}{\varphi}(F_{2n+1}+F_{2n+3})\right]$\\

Since we know, $F_{2n}+F_{2n+2}=L_{2n+1}$ and $F_{2n+1}+F_{2n+3}=L_{2n+2}$, where $L_n$ are Lucas numbers;\\

So,$$A_{n,n+1}+A_{n+1,n+2}=\frac{\pi}{20ln(\varphi)}\left[L_{2n+1}+\varphi'L_{2n+2}\right].$$\\
Since $L_{n}=\varphi^{n}+\varphi'^{n},$
\begin{align}
L_{2n+1}+\varphi'L_{2n+2}&=\varphi^{2n+1}+\varphi'^{2n+1}+\varphi'(\varphi^{2n+2}+\varphi'^{2n+2})\nonumber \\
&=\varphi^{2n+1}+\varphi'^{2n+1}+\varphi'\varphi\varphi^{2n+1}+\varphi'^{2n+3} \;\;\mbox{since}\,\, \varphi'\varphi=-1 ;\nonumber \\
&=\varphi'^{2n+1}+\varphi'^{2n+3}\nonumber \\
&=\varphi'^{2n+1}(1+\varphi'^2) \;\;\mbox{and since}\,\,\varphi'^2=\varphi'+1\nonumber \\
&=\varphi'^{2n+1}(2+\varphi')\nonumber \\
\nonumber
\end{align}
Thus, we obtained;\\
$L_{2n+1}+\varphi'L_{2n+2}=\varphi^{2n+1}(2+\varphi')=\left(\frac{-1}{\varphi}\right)^{2n+1}\left(2-\frac{1}{\varphi}\right)=\fbox{$-(\frac{1}{\varphi})^{2n+1}(2-\frac{1}{\varphi})$}$\\

Now, for both\,\, $n=2k$ and $n=2k+1$;\quad $-\left(\frac{1}{\varphi}\right)^{2n+1}\left(2-\frac{1}{\varphi}\right)<0$.

Thus;\\

$A_{n,n+1}+A_{n+1,n+2}=\frac{\pi}{20ln(\varphi)}\left[L_{2n+1}+\varphi'L_{2n+2}\right]<0$ \\

 The proof is completed.\\
\phantom{abcdefghklmoprstuyzbmerveozvatanmerveozvatanmervemerve}$\blacksquare$\\

\subsection{Relative Area Sequence}

Here, we define $\lambda_n$ as ratio of areas;\\

$$\lambda_n=\left|\frac{A_{n,n+1}}{A_{n-1,n}}\right|=\left|\frac{8(ln(\varphi))^2(-1)^n-\pi^2(F_{2n}-\frac{1}{\varphi}F_{2n+1})}{-8(ln(\varphi))^2(-1)^n-\pi^2(F_{2n-2}-\frac{1}{\varphi}F_{2n-1})}\right|$$\\

$$=\frac{|(F_{2n}-\frac{1}{\varphi}F_{2n+1})-\frac{8(ln(\varphi))^2}{\pi^2}(-1)^n|}{|(F_{2n-2}-\frac{1}{\varphi}F_{2n-1})+\frac{8(ln(\varphi))^2}{\pi^2}(-1)^n|}$$\\

$$=\frac{|(\varphi')^{2n+1}-\frac{8(ln(\varphi))^2}{\pi^2}(-1)^n|}{|(\varphi')^{2n-1}+\frac{8(ln(\varphi))^2}{\pi^2}(-1)^n|}$$\\

$$=\frac{|(-1)^{2n+1}\frac{1}{\varphi^{2n+1}}-\frac{8(ln(\varphi))^2}{\pi^2}(-1)^n|}{|(-1)^{2n-1}\frac{1}{\varphi^{2n-1}}+\frac{8(ln(\varphi))^2}{\pi^2}(-1)^n|}$$\\

$$=\frac{|-\frac{1}{\varphi^{2n+1}}-\frac{8(ln(\varphi))^2}{\pi^2}(-1)^n|}{|-\frac{1}{\varphi^{2n-1}}+\frac{8(ln(\varphi))^2}{\pi^2}(-1)^n|}$$\\

$$=\frac{|\frac{1}{\varphi^{2n+1}}+\frac{8(ln(\varphi))^2}{\pi^2}(-1)^n|}{|\frac{1}{\varphi^{2n-1}}-\frac{8(ln(\varphi))^2}{\pi^2}(-1)^n|}$$\\

$$=\frac{\varphi^{2n-1}}{\varphi^{2n+1}}\frac{|1+\frac{8(ln(\varphi))^2}{\pi^2}(-1)^n\varphi^{2n+1}|}{|1-\frac{8(ln(\varphi))^2}{\pi^2}(-1)^n\varphi^{2n-1}|}$$\\

$$=\frac{1}{\varphi^2}\frac{|\frac{8(ln(\varphi))^2}{\pi^2}\varphi^{2n}\varphi+(-1)^n|}{|\frac{8(ln(\varphi))^2}{\pi^2}\frac{\varphi^{2n}}{\varphi}-(-1)^n|}$$\\

$$\mbox{At the end,}\;\lambda_n=\left|\frac{A_{n,n+1}}{A_{n-1,n}}\right|=\frac{1}{\varphi^2}\frac{|\frac{8(ln(\varphi))^2}{\pi^2}\varphi^{2n}\varphi+(-1)^n|}{|\frac{8(ln(\varphi))^2}{\pi^2}\varphi^{2n}\frac{1}{\varphi}-(-1)^n|}\;\;\mbox{is obtained.} $$\\

\end{document}